\documentclass[11pt]{article}

\usepackage{color}
\usepackage{amsmath}
\usepackage{amssymb}
\usepackage{indentfirst}
\usepackage{graphics} 
\usepackage{color}
\usepackage{cases}
\makeatletter

\@addtoreset{equation}{section}
\makeatother
\def\<{\langle}
\def\>{\rangle}

\newtheorem{lem}{Lemma}[section]
\newtheorem{theo}{Theorem}[section]
\newtheorem{rem}{Remark}[section]
\newtheorem{pro}{Proposition}[section]

\newtheorem{cor}{Corollary}[section]
\makeatletter
   
   \@addtoreset{equation}{section}
\makeatother

\begin{document}
\title{\bf Optimal $L^{2}$-growth of the generalized\\ Rosenau equation}
\author{Xiaoyan Li\thanks{xiaoyanli@hust.edu.cn} \\ {\small School of Mathematics and Statistics}\\ {\small Huazhong University of Science and Technology} \\ {\small Wuhan, Hubei 430074, PR China}\\ and\\Ryo Ikehata\thanks{Corresponding author:ikehatar@hiroshima-u.ac.jp} \\ {\small Department of Mathematics, Division of Educational Sciences}\\ {\small Graduate School of Humanities and Social Sciences} \\ {\small Hiroshima University} \\ {\small Higashi-Hiroshima 739-8524, Japan}}
\date{}

\date{}
\maketitle
\begin{abstract}
We report that the quantity measured in the $L^2$ norm of the solution itself of the generalized Rosenau equation, which was completely unknown in this equation, grows in the proper order at time infinity. It is also immediately apparent that this growth aspect does not occur in three or more spatial dimensions, so we will apply the results obtained in this study to provide another proof that Hardy-type inequalities do not hold in the case of one or two spatial dimensions.
\end{abstract}
\section{Introduction}
\footnote[0]{Keywords and Phrases: Rosenau equation; Fourier analysis; weighted $L^{1}$-data; low dimensional case; growth estimates.}
\footnote[0]{2020 Mathematics Subject Classification. Primary 35B40, 35B45; Secondary 35B05, 35E15 35G10, 35R11.}

We consider the Cauchy problem of of the generalized Rosenau equation: 
\begin{numcases}
  ~u_{tt} +\delta (- \Delta)^\theta u_{tt} +\mu  \Delta^2u-\kappa  \Delta u = 0,\ \ \ (t,x)\in (0,\infty)\times {\bf R}^{n},\label{d1}\\
 u(0,x)= u_0(x), \quad  u_{t}(0,x)= u_{1}(x),\ x\in{\bf R}^{n}, \label{d2}
\end{numcases}
where $\delta >0$,~$\mu>0$, $\kappa>0$, and $0 < \theta \leq 2$. Incidentally, we call it so because our target is the case $\theta = 2$. We assume that, for the moment, $[u_0(x), u_1(x) ]\in H^{2}({\bf R}^{n})\times H^{\theta}({\bf R}^{n})$. 

The total energy $E(t)$ for the solution to problem \eqref{d1}-\eqref{d2} is defined by
\begin{equation}\label{d3}
		E(t) := \frac{1}{2}\left(\Vert u_{t}(t,\cdot)\Vert_{L^{2}({\bf R}^{n})}^{2} +\delta  \Vert (-\Delta)^{\frac{\theta}{2}} u_{t}(t,\cdot)\Vert_{L^{2}({\bf R}^{n})}^{2} 
		+\mu  \Vert\Delta u(t,\cdot)\Vert_{L^{2}({\bf R}^{n})}^{2}
		+\kappa  \Vert\nabla u(t,\cdot)\Vert_{L^{2}({\bf R}^{n})}^{2}\right).
\end{equation}
Taking derivative of $t$ on the both sides of \eqref{d3}, we see that 
\[
\frac{d }{dt}E(t)=0,
\]
which implies the solution of \eqref{d1}-\eqref{d2} satisfies energy conservation law $E(t)=E(0)$ for $t\geq 0$. However, we are very interested in the time behavior of the solution itself, since energy itself does not directly contain information about the solution itself.

Let us start with the Rosenau equation proposed by P. Rosenau \cite{rp1968, rp1968-2}, which is one of main topics for mathematicians and physicists. The  global existence and uniqueness of  solution to the first order in time Rosenau equation were demonstrated by \cite{p1999} (see also \cite{p1999-2}).  Concerning the second order in time Rosenau equation
\begin{align} \label{d61}
	u_{tt}+u_{xxxx}+u_{xxxxtt}-\gamma u_{xx}=f(u)_{xx},
\end{align}
the authors in \cite{wangxu09} proved the existence and uniqueness of global solution under the assumption $f(u)=-\beta |u|^pu$ ($\beta>0,~ p>0$), and discussed the finite-time blow-up of solution with the aid of potential well method. When the hydrodynamical friction was taken into consideration, more general Rosenau equation 
\begin{align} \label{d60}
	u_{tt}-\Delta u-\nu \Delta  u_{t}+ \Delta^2 u+ \Delta^2 u_{tt}= \Delta f(u)
\end{align}
appeares. For \eqref{d60}, the global
existence and asymptotic behavior of the solution were showed  in the time-weighted Sobolev spaces under smallness condition on the initial data \cite{ws2013}. Note that \eqref{d60} and \eqref{d61} are both nonlinear problems. Originally,  Rosenau equations are devoted to some high
order nonlinear systems \cite{rp1968}.
However, it is natural that the results of  corresponding linear part (i.e. $f(u)=0$) have some influence on the behavior of solutions to nonlinear problem because of Duhamel’s principle, such as the analyses in \cite{ws2013}.
Therefore, it is worth to consider the linear Rosenau equation (i.e. $f(u)=0$) and there are some interesting consequences related to linear Rosenau equation. For example, the author in \cite{mh2019} considers the initial value problem for 
\begin{align} \label{d63}
	u_{tt}-\Delta u-\nu \Delta u_t+\Delta^2 u+\Delta^2 u_{tt}=0,
\end{align}
and presented the asymptotic estimates of solutions. Note that the equation \eqref{d63} has a dissipative structure. In some literatures, \eqref{d63} is also called the plate equation with rotational structure inertia (see \cite{H2016}) or generalized Boussinesq equation with
double dispersive term (see \cite{pn2009}).
Such higher order hyperbolic equations with some dissipative structures are also studied by
many mathematicians (see \cite{cr2013, HIC, cr2009, ir2015} for plate equations and \cite{th2012, lz1999} for beam equations).

In this paper, we consider the Cauchy problem for the linear Rosenau equation \eqref{d1}
and establish the optimal $L^2$-growth property of the solution itself in the case of $n = 1,2$. To the best of our knowledge, to capture the optimal $L^2$-growth estimate of the solution itself in \eqref{d1} is still an open problem, and it shows a singularity of the solution itself near low frequency zone. The Rosenau equation without any dissipative terms is a class of important partial differential equations, whatever the case, it should be studied.

The motivation of  this paper includes two  aspects. On the one hand, some traditional $L^p$-$L^q$ estimates are a bit far from the best estimate when one tries to get the $L^2$-estimate of the solution itself, such as the estimates of free waves in \cite{B, M, P, Pl, S, W} and the references therein. When strong damping are added into free waves, the results about best growth estimates of various wave-derived equations are obtained (see \cite{Ba, CI, CT, FIM, I-14, IT, IO, Mi}), and Ikehata's optimal estimates of the  $L^2$ norm for free wave \cite{JHDE-ike} and plate equations \cite{FE-ike} are one of the most fundamental results in these studies. 

On the other hand, in a series of papers (see, for example, \cite{rcc2020, ike2005}), the core idea to get the boundedness of $L^2$-norm of the solution itself is to use the Hardy type inequality as will be discussed in Section 5 below. However, in this paper, 
we analyze the $L^2$-estimate of the solution in low-dimensional case ($n = 1,2$) in ${\bf R}^{n}$. In this case, the Hardy type inequality (for $n=2$) never holds as was already pointed out in (at least) \cite{MOW}, so it requires more delicate treatments to treat our problem.

 Before going to introduce our theorems, we set the following notations.
\\

{\bf Notation.} {\small Throughout this paper, $\| \cdot\|_q = \Vert\cdot\Vert_{L^{q}({\bf R}^{n})}$ stands for the usual $L^q({\bf R}^{n})$-norm. For simplicity of notation, in particular, we use $\| \cdot\|$ instead of $\| \cdot\|_2$. Furthermore, $\Vert\cdot\Vert_{H^{s}({\bf R}^{n})}$ stands for the usual $H^{s}$-norms of the Sobolev spaces $H^{s}({\bf R}^{n})$ ($s \geq 0$). We also introduce the following weighted functional spaces
	\[L^{1,\gamma}({\bf R}^{n}) := \left\{f \in L^{1}({\bf R}^{n}) \;\big| \; \Vert f\Vert_{1,\gamma} := \int_{{\bf R}^{n}}(1+\vert x\vert^{\gamma})\vert f(x)\vert dx < +\infty\right\}.\]
The Fourier transform ${\cal F}_{x\to\xi}(f)(\xi)$ of $f(x)$ is defined by 
	\[{\cal F}_{x\to\xi}(f)(\xi) = \hat{f}(\xi) := \displaystyle{\int_{{\bf R}^{n}}}e^{-ix\cdot\xi}f(x)dx, \quad \xi \in {\bf R}^n,\]
	as usual with $i := \sqrt{-1}$, and ${\cal F}_{\xi\to x}^{-1}$ expresses its inverse Fourier transform. The fractional operators $(-\Delta)^{\theta}$ ($\theta \geq 0$) can be defined by (formally)
\[\left((-\Delta)^{\theta}u\right)(x) := {\cal F}_{\xi\to x}^{-1}(\vert\xi\vert^{2\theta}\hat{u}(\xi))(x).\]
We denote the surface area of the $n$-dimensional unit ball by $\omega_{n} := \displaystyle{\int_{\vert\omega\vert = 1}}d\omega$.
For each $n = 1,2$, we set
\[I_{0,n} := \Vert u_{1}\Vert_{L^{2}({\bf R}^{n})} + \Vert u_{1}\Vert_{L^{1}({\bf R}^{n})}.\]
The inner product of $L^2({\bf R}^n)$ is defined by
$$(f,g)=\int_{{\bf R}^{n}} f(x)\overline{g(x)} dx,~~~f,g\in L^2({\bf R}^n).$$ 
The equivalent inner products of  $H^\theta({\bf R}^{n}) \: (\theta>0) $ are denoted by 
\[ (u,v)_{H^\theta}=\int_{{\bf R}^n_\xi} (1+ \delta \vert\xi\vert^{2\theta})
\hat{u}\bar{\hat{v}} ~d\xi  \cong (u,v)+\delta \big((-\Delta)^{\frac{\theta}{2}}u, (-\Delta)^{\frac{\theta}{2}}v\big),
\]
\[
(u,v)_{H^2}=\int_{{\bf R}^n_\xi} (1+ \kappa |\xi|^{2}+\mu|\xi|^4)
\hat{u}\bar{\hat{v}} ~d\xi  \cong (u,v)+\kappa \big((-\Delta)^{\frac{1}{2}}u, (-\Delta)^{\frac{1}{2}}v \big)+\mu\big((-\Delta)u, (-\Delta)v \big).
\]
}
At any rate we proceed on the basis of the following proposition concerning the existence and uniqueness of solutions. We refer the reader to Appendix for an outline of the proof. The proof itself as stated in Appendix is just a simple modification of \cite{H2016, mh2019}.
\begin{pro}\label{prp1}
Let $n \geq 1$ and $\theta \in [0,2]$. If $[u_0, u_1]\in H^2({\bf R}^n)\times H^\theta({\bf R}^n) $, the equation \eqref{d1}-\eqref{d2} admits a unique mild solution 
	$$u \in  C([0,\infty); H^{2}({\bf R}^n)) \cap C^1([0,\infty); H^{\theta}({\bf R}^n)).$$
If the initial data further satisfies $[u_0, u_1]\in H^{4-\theta}({\bf R}^n)\times H^2({\bf R}^n) $, the equation \eqref{d1}-\eqref{d2} admits a unique strong solution 
$$u \in  C([0,\infty); H^{4-\theta}({\bf R}^n)) \cap C^1([0,\infty); H^{2}({\bf R}^n))\cap C^2([0,\infty); H^{\theta}({\bf R}^n)).$$
\end{pro}

Our first result is concerned with the optimal growth property in the case of $n = 1$.

\begin{theo}\label{theorem1}
Let $n = 1$ and $0 < \theta \leq 2$, $u_{0} \in H^2({\bf R})$ and $u_{1} \in H^{\theta}({\bf R})$. Then, the solution $u(t,x)$ to problem \eqref{d1}-\eqref{d2} satisfies the following properties under the additional regularity on the initial data:
\[ \Vert u(t,\cdot)\Vert_{L^{2}({\bf R})} \leq C_{1}I_{0,1}\sqrt{t},~~~~~\text{if}~~
u_{1} \in L^{1}({\bf R}),\]
\[ C_{2}\left\vert\int_{{\bf R}}u_{1}(x)dx\right\vert\sqrt{t} \leq \Vert u(t,\cdot)\Vert_{L^{2}({\bf R})}, ~~~~~\text{if}~~
u_{1} \in L^{1,\gamma}({\bf R}),~\gamma \in (\frac{1}{2},1],\]
for $t \gg 1$, where $C_{j} > 0$ {\rm ($j = 1,2$)} are constants depending only on the space dimension and $\gamma$.
\end{theo}
Our next result is the case of $n = 2$.
\begin{theo}\label{theorem2}
Let $n = 2$ and $0 < \theta \leq 2$,  $u_{0} \in H^2({\bf R}^{2})$ and $u_{1} \in H^{\theta}({\bf R}^{2})$. Then, the solution $u(t,x)$ to problem \eqref{d1}-\eqref{d2} satisfies the following properties under the additional regularity on the initial data:
\[ 
\Vert u(t,\cdot)\Vert_{L^{2}({\bf R}^{2})} \leq C_{1}I_{0,2}\sqrt{\log t}, ~~~\text{if}~~u_{1} \in L^{1}({\bf R}^{2}),\]
\[
 C_{2}\left\vert\int_{{\bf R}^{2}}u_{1}(x)dx\right\vert\sqrt{\log t} \leq \Vert u(t,\cdot)\Vert_{L^{2}({\bf R}^{2})},     ~~~~~\text{if}~~
 u_{1} \in L^{1,\gamma}({\bf R}^{2}),~\gamma \in (0,1],
  \]
for $t \gg 1$, where $C_{j} > 0$ {\rm ($j = 1,2$)} are constants depending only on the space dimension and $\gamma$. 
\end{theo}
\begin{rem}{\rm Comparing the above two theorems with the previous study \cite{JHDE-ike, FE-ike}, it can be seen that when viewed in terms of the growth of the $L^2$-norm of the solution at time infinity that the equation has, it originates from a wave equation with $\mu = 0$ and $\delta = 0$, and the influence of the Plate equation with $\delta = \kappa = 0$ disappears. The final growth property is due to the influence from the low-frequency band of the solution. The "dissipative structure" of the solution due to the influence from the low-frequency band can be seen in the paper by \cite{th2012}, who studied equation \eqref{d1} with $\delta = 0$ and additional friction term.
Incidentally, the results of \cite{IL} are helpful for the growth estimates when $\mu=0$, $\kappa = 1$, $\delta = 1$ and $\theta = 1$.}
\end{rem}
\begin{rem}{\rm The assumption of a class of initial functions $u_{0} \in H^2({\bf R}^{n})$ and $u_{1} \in H^{\theta}({\bf R}^{n})$ ($n = 1,2$) is sufficient information for the unique existence of a mild solution as can be seen from Proposition \ref{prp1} and is not directly involved in the growth estimate itself. }
\end{rem}

This paper is organized as follows. The proofs of the estimate formulas from below for Theorems \ref{theorem1} and \ref{theorem2} will be discussed in Section 2, and the proofs of the estimates from above will be derived in Section 3. We note in Section 4 that $L^2$-growth does not occur in three or more dimensions. In Section 5, as a remark, we discuss the validity of the Hardy type inequalities in the low dimensional case. Appendix describes the outline of Proposition \ref{prp1}.



\section{Estimates from below of the solution: $n = 1,2$}

In this section, we derive the $L^2$-lower bound of solution for $n = 1,2$. By Plancherel theorem, we just need to discuss the  the $L^2$-lower bound of $w(t,\xi)$ defined by \eqref{d36}. The main task is to deal with the estimates in the low frequency domain ($|\xi|\leq \beta (t)$), where $\beta(t)$ is a time-variable function chosen suitably.\\
Before those proofs, prepare the following three tools and facts required in this section.
Define (possibly $L = 1$)
\begin{equation}\label{d50}
L := \sup_{\eta \ne 0}\left\vert \frac{\sin\eta}{\eta}\right\vert < +\infty.
\end{equation}
Let $\delta_{0} \in (0,1)$ be a real number such that the following inequality holds:
\begin{equation}\label{d31}
\left\vert \frac{\sin\eta}{\eta}\right\vert \geq \frac{1}{2}
\end{equation}
for all $\eta \in (0,\delta_{0}]$. We finally prepare the inequality
\begin{equation}\label{d33}
\vert a + b\vert^{2} \geq \frac{1}{2}\vert a\vert^{2} - \vert b\vert^{2} 
\end{equation}
for all $a, b \in {\bf C}$.\\

In the Fourier space ${\bf R}_{\xi}^{n}$ the problem \eqref{d1}-\eqref{d2} can be transformed into the following ODE with parameter $\xi \in {\bf R}_{\xi}^{n}$
\begin{numcases}
 ~(1+\delta \vert\xi\vert^{2\theta})w_{tt} +(\mu |\xi|^4+ \kappa |\xi|^2)w = 0,\ \ \ t>0,\quad \xi \in {\bf R}_{\xi}^{n},\label{5}\\
 w(0,\xi) = w_0(\xi), \quad  w_{t}(0,\xi)= w_{1}(\xi),\ \ \ \xi \in{\bf R}^{n} ,\label{6}
\end{numcases}
where  $w_{0}(\xi) := \hat{u}_0(\xi)$, $w_{1}(\xi) := \hat{u}_1(\xi)$ and $w(t,\xi) := \hat{u}(t,\xi)$. Moreover, the solution to  problem \eqref{5}-\eqref{6} is written explicitly in the form of 
\begin{equation}\label{d36}
w(t,\xi) =\cos(t f(\vert\xi\vert)  ) w_{0}(\xi)+ \frac{\sin(tf(\vert\xi\vert))}{f(\vert\xi\vert)}w_{1}(\xi),
\end{equation}
where we set
\begin{equation}\label{d55}
f(r) := \sqrt{\frac{\mu r^4+\kappa r^2}{1+\delta r^{2\theta}}}.
\end{equation}
With these preparations in place, the proofs will begin immediately in this section and in Section 3.\\

In this section, we will consider the $1$-D case, but for the time being, we will proceed with the discussion in $n$-dimensions and limit it to $1$-D when it becomes necessary. By density, we can assume that the initial data $u_{0}, u_{1} \in C_{0}^{\infty}({\bf R}^{n})$.
Then we define a subset $L_{0}$ of ${\bf R}_{\xi}^{n}$ with $n\geq 1$ as follows:
\begin{equation}\label{d30}
L_{0} := \left\lbrace \xi \in {\bf R}_{\xi}^{n}\,:\,\vert\xi\vert \leq \beta(t) := \frac{\delta_{0}}{(\mu+\kappa)^{\frac{1}{2}}  t} \right\rbrace ,
\end{equation}
with $\delta_0$ defined in \eqref{d31}, and take $t\geq \frac{\delta_0}{\sqrt{\mu+\kappa}}$ sufficiently large 
guaranteeing $\beta(t) \leq 1$.
By the definition of $f(r)$ in \eqref{d55}, we have that
\begin{align} \label{d70}
\xi \in L_0 \quad \Rightarrow \quad
\big(\frac{\delta_0}{t}\big)^2 \geq (\mu+\kappa)|\xi|^2\geq 
\frac{\mu |\xi|^4+\kappa|\xi|^2}{1+\delta |\xi|^{2\theta}} \quad \Rightarrow \quad tf(\vert\xi\vert) \in [0,\delta_{0}],~~~~(|\xi|\leq 1).
\end{align}
Using \eqref{d31} and \eqref{d70}, one can estimate the solution $w(t,\xi)$ in the case for $n=1$ as follows:
\begin{equation}\label{d35}
I_{l}(t) := \int_{L_{0}}\frac{\sin^{2}(tf(\vert\xi\vert))}{f^{2}(\vert\xi\vert)}d\xi \geq \frac{t^{2}}{4}\int_{L_{0}}d\xi  = \frac{t^{2}}{2}\frac{\delta_{0}}{(\mu+\kappa)^{\frac{1}{2}}  t} \geq C t
\end{equation}
for $t \gg 1$, where $C > 0$ is some constant. \\

Next we decompose the initial data $w_{1}(\xi)$ in the Fourier space
\begin{equation}\label{d34}
w_{1}(\xi) = P + (A(\xi)-iB(\xi)),\quad \xi \in {\bf R}_{\xi}^{n},\quad (n \geq 1)
\end{equation}
where
\[P := \int_{{\bf R}^{n}}u_{1}(x)dx,\]
\[A(\xi) := \int_{{\bf R}^{n}}(\cos(x\xi)-1)u_{1}(x)dx, \quad B(\xi) := \int_{{\bf R}^{n}}\sin(x\xi)u_{1}(x)dx.\]
Using the results in \cite{I-04}, there is a constant $M > 0$ such that 
\begin{equation}\label{d32}
\vert A(\xi)-iB(\xi)\vert \leq M\vert\xi\vert^{\gamma}\Vert u_{1}\Vert_{1,\gamma},\quad \xi \in {\bf R}_{\xi}^{n},
\end{equation}
where $u_{1} \in L^{1,\gamma}({\bf R}^{n})$ and $\gamma \in (0,1]$. Combinig \eqref{d33} with \eqref{d34} with $n = 1$, we get
\begin{align} \label{d37}
	J_{1}(t)&: = \frac{1}{2}\int_{L_{0}}\frac{\sin^{2}(tf(\vert\xi\vert))}{f^2(\vert\xi\vert)}\vert w_{1}(\xi)\vert^{2}d\xi-\int_{{\bf R}} \left| \cos(t f(\vert\xi\vert) ) \right| ^2 \left| w_{0}(\xi)\right|^2 ~d\xi  \notag \\
	&\geq \frac{P^{2}}{4}\int_{L_{0}}\frac{\sin^{2}(tf(\vert\xi\vert))}{f^2(\vert\xi\vert)}d\xi - \int_{L_{0}}\vert A(\xi)-iB(\xi)\vert^{2}\frac{\sin^{2}(tf(\vert\xi\vert))}{f^2(\vert\xi\vert)}d\xi - \int_{{\bf R}}  \left| w_{0}(\xi)\right|^2 ~d\xi \notag \\
	&=: \frac{P^{2}}{4}I_{l}(t)-R_{l}(t)-\|w_0(\xi)\|^2.
\end{align}
By \eqref{d32},  for $\gamma \in (1/2,1]$,  $R_{l}(t)$ can be  estimated by
\begin{align}	\label{d38}
	R_{l}(t) &\leq M^{2}\Vert u_{1}\Vert_{1,\gamma}^{2}\int_{L_{0}}\frac{r^{2\gamma}}{f^{2}(r)}~d\xi\notag  \\
	&= M^{2}\Vert u_{1}\Vert_{1,\gamma}^{2}\int _{L_0} \frac{(1+\delta r^{2\theta })r^{2\gamma}}{\mu r^4+\kappa r^2 }~dr \notag \\
		&\leq 2\frac{1+\delta}{\kappa } M^2\Vert u_{1}\Vert_{1,\gamma}^{2} \int _0^{\beta(t)} r^{2\gamma-2}~dr\notag\\
			&\leq \frac{1+\delta}{\kappa(2\gamma-1)  } M^2\Vert u_{1}\Vert_{1,\gamma}^{2} \left[r^{2\gamma-1}\right]_0^{\beta(t)}\notag \\
				&\leq C_\gamma  \Vert u_{1}\Vert_{1,\gamma}^{2} ~t^{-(2\gamma-1)},
\end{align}
where $r=|\xi|$,   $C_\gamma>0$ and $t \gg 1$, and $2\gamma > 1$. Thus, from \eqref{d36}, \eqref{d33}, \eqref{d35}, \eqref{d37} and \eqref{d38}, we can get the desired lower bound estimate
\[\Vert w(t,\cdot)\Vert^{2} \geq J_{1}(t) \geq CP^{2}t - C_{\gamma}\Vert u_{1}\Vert_{1,\gamma}^{2}t^{-(2\gamma-1)}-\|u_0\|^2 , \quad t \gg 1\]
with some constants $C > 0$ and $C_{\gamma} > 0$. By density, we have the following lemma.
\begin{lem}\label{lem1}Let $n = 1$, $\gamma \in (\frac{1}{2},1]$, and $\mu \geq 0$, $\theta \in [0,2]$ and $\kappa > 0$. Assume $u_{0} \in L^2({\bf R})$ and $u_{1} \in L^{1,\gamma}({\bf R})$. Then, it holds that
\[\Vert w(t,\cdot)\Vert^{2} \geq CP^{2} t,
 \quad~~ (t \gg 1).\]
\end{lem}


Next, we discuss the  two dimensional case.
In this part, by  integration by parts, we make full use of the characteristics of $f(r)$ to get the desired estimate of solution.

We start with the integral combining the trick function $e^{-\vert\xi\vert^{2}}$. This idea comes from the paper \cite{JHDE-ike}.  

It follows from  \eqref{d33}, \eqref{d34} and \eqref{d32}  that 
\begin{align} 
	\Vert w(t,\cdot)\Vert^{2} &\geq \frac{1}{2}\int_{{\bf R}^{2}}\frac{\sin^{2}(tf(\vert\xi\vert))}{f^{2}(\vert\xi\vert)}\vert w_{1}(\xi)\vert^{2}d\xi-
	 \int_{{\bf R}^{2}} |\cos^{2}(tf(\vert\xi\vert))|^2 | w_{0}(\xi)|^2 d\xi
	\notag\\
	&\geq \frac{1}{2}\int_{{\bf R}^{2}}e^{-\vert\xi\vert^{2}}\frac{\sin^{2}(tf(\vert\xi\vert))}{f^{2}(\vert\xi\vert)}\vert P+(A(\xi)-iB(\xi))\vert^{2}d\xi
	-
	\int_{{\bf R}^{2}} | w_{0}(\xi)|^2 d\xi
	\notag\\
	&\geq \frac{1}{4}P^{2}\int_{{\bf R}^{2}}e^{-\vert\xi\vert^{2}}\frac{\sin^{2}(tf(\vert\xi\vert))}{f^{2}(\vert\xi\vert)}d\xi - \int_{{\bf R}^{2}}e^{-\vert\xi\vert^{2}}\frac{\sin^{2}(tf(\vert\xi\vert))}{f^{2}(\vert\xi\vert)}\left(M^{2}\Vert u_{1}\Vert_{1,\gamma}^{2}\vert\xi\vert^{2\gamma}\right)d\xi\notag-\|w_0\|^2\\
	&\geq \frac{1}{4}P^{2}\int_{{\bf R}^{2}}e^{-\vert\xi\vert^{2}}\frac{\sin^{2}(tf(\vert\xi\vert))}{f^{2}(\vert\xi\vert)}d\xi - M^{2}\Vert u_{1}\Vert_{1,\gamma}^{2}\int_{{\bf R}^{2}}e^{-\vert\xi\vert^{2}}\frac{\vert\xi\vert^{2\gamma}}{f^{2}(\vert\xi\vert)}d\xi-\|u_0\|^2\notag\\
	&=:  \frac{1}{4}P^{2}T(t) -M^{2}\Vert u_{1}\Vert_{1,\gamma}^{2}U(t)-\|u_0\|^2.\label{d71}
\end{align}
The estimate of $U(t)$ are given by
\begin{align}\label{i-1}
	\frac{U(t)}{\omega_{2}}
	&= \int_{0}^{\infty}e^{-r^{2}}r^{2\gamma+1}\frac{1+\delta r^{2\theta}}{\mu r^4+\kappa r^2}dr \nonumber\\
	&\leq\frac{1}{\kappa} \int_{0}^{\infty}e^{-r^{2}}r^{2\gamma-1}dr + \frac{\delta}{\kappa} \int_{0}^{\infty}e^{-r^{2}}r^{2\gamma+2\theta-1}dr\nonumber\\
	&= \frac{1}{\kappa\gamma}\left( \int_{0}^{\infty}e^{-r^{2}}r^{2\gamma+1}dr +  \left[ \frac{1}{2}
	r^{2\gamma}e^{-r^2}\right]_0^\infty\right) \nonumber\\
     &+\frac{1}{\kappa(\gamma+\theta)} \left( \int_{0}^{\infty}
	e^{-r^{2}}r^{2(\gamma+\theta)+1}dr+
	 \left[ \frac{1}{2}
	r^{2(\gamma+\theta)}e^{-r^2}\right]_0^\infty\right) \nonumber\\
	&= \frac{1}{\kappa\gamma} \int_{0}^{\infty}e^{-r^{2}}r^{2\gamma+1}dr+
		\frac{1}{\kappa(\gamma+\theta)}  \int_{0}^{\infty}
	e^{-r^{2}}r^{2(\gamma+\theta)+1}dr\nonumber\\
	&=: K_{0}.
\end{align}
Note that for $\gamma \in (0,1]$, $K_{0}$ is finite valued . Here the trick function $e^{-\vert\xi\vert^{2}}$ plays an important role.

Before providing  the estimate of $T(t)$, the characteristics of $f(r)$ are shown by the following proposition, which is useful to get more detailed estimates of $T(t)$.

\begin{pro} \label{d40}
Let $\delta  > 0$, $2 \geq \theta > 0$, $\mu \geq 0$ and $\kappa > 0$. Setting
\[\varepsilon_0=\min \left\lbrace \left( \frac{\kappa}{2(\mu+\kappa)\delta \theta}\right)  ^\frac{1}{2\theta},~1\right\rbrace, 
\]	
for $0 < \forall r\leq \varepsilon_0\leq1$, it holds that
\[
\left| \frac{1}{f'(r)}\right| \leq C ~~~~~\text{and} ~~~~~\vert f''(r)\vert \leq C (1+\frac{1}{r}),
\]	
where $C > 0$ is an universl constant depending on the parameters $\delta  > 0$, $\theta > 0$, $\mu \geq 0$ and $\kappa > 0$.
\end{pro}
\noindent{ \it{Proof.}}
It follows from the definition of $f(r)$ in \eqref{d55} that

\begin{align}
\left| f^{'}(r)\right| &=\frac{	\left| (4\mu r^{3}+2\kappa r)(1+\delta r^{2\theta})-2\delta  \theta(\mu r^{4}+\kappa r^{2}) r^{2\theta-1}\right| }{
	2(\mu r^{4}+\kappa r^{2})^\frac{1}{2}(1+\delta r^{2\theta})^{\frac{3}{2}}	
}\notag\\
	&\geq \frac{	\left| (4\mu r^{3}+2\kappa r)(1+\delta r^{2\theta})-2\delta  \theta(\mu r^{4}+\kappa r^{2}) r^{2\theta-1}\right| }{
		2(\mu+\kappa)^\frac{1}{2} (1+\delta)^\frac{3}{2} r	
	}\notag \\
	&\geq \frac{ 2\kappa r-2\delta  \theta(\mu +\kappa ) r^{2\theta+1} }{
		2(\mu+\kappa)^\frac{1}{2} (1+\delta)^\frac{3}{2} r	
	}\notag \\
&=\frac{ \kappa -\delta  \theta(\mu +\kappa ) r^{2\theta} }{
	(\mu+\kappa)^\frac{1}{2} (1+\delta)^\frac{3}{2} 	
}\notag \\
	&>\frac{\kappa}{2} \frac{1}{(\mu+\kappa)^\frac{1}{2} (1+\delta)^\frac{3}{2}},  \quad (0< \forall r \leq \varepsilon_{0} \leq 1),
\end{align}
which implies the boundedness of $\left|  \frac{1}{f^{'}(r)}\right| $.


By tedious but simple calculations, for  $0 < r\leq \varepsilon_0\leq1$, $f''(r)$ satisfies
\begin{align}
	\left| f''(r)\right| =& \bigg|	\frac{\left(6 \mu r^2+
		\kappa \right)\left(\mu r^4+\kappa r^2\right)^{\frac{1}{2}}\left(1+\delta r^{2 \theta}\right)^{\frac{1}{2}} }{\left(\mu r^4+\kappa r^2\right)\left(1+\delta r^{2 \theta}\right)} \notag \\
	&-	\frac{\left(2 \mu r^3+\kappa r\right)\left(4 \mu r^3+2 \kappa  r\right)\left(1+\delta r^{2 \theta}\right)+\left(2 \mu r^3+\kappa r\right)\left(\mu r^4+\kappa r^2\right)\left(2 \theta \delta r^{2 \theta-1}\right)}{2\left[\left(\mu r^4+\kappa r^2\right)\left(1+\delta r^{2 \theta}\right)\right]^{\frac{3}{2}} }\notag\\
	&	- \frac{\delta \theta \left(4 \mu r^3+2\kappa  r\right) r^{2 \theta-1} }{2 \left(\mu r^4+\kappa r^2\right)^{\frac{1}{2}}\left(1+\delta r^{2 \theta}\right)^{\frac{3}{2}} }-\frac{\delta \theta(2 \theta-1) \left(\mu r^4+\kappa r^2\right)^{\frac{1}{2}}  r^{2 \theta-2} }{\left(1+\delta \theta r^2 \right)^{\frac{3}{2}}}+\frac{ 3\delta^2 \theta^2  \left(\mu r^4+\kappa r^2\right)^{\frac{1}{2}}   r^{4 \theta-2} }{\left(1+\delta r^{2 \theta}\right)^{\frac{5}{2}}} \bigg|\notag\\
	\leq& 
	\frac{(6\mu+\kappa )(\mu+\kappa )^{\frac{1}{2}} (1+\delta)^{\frac{1}{2}}r}{\kappa r^2}+\frac{(2\mu+\kappa )(4\mu+2 \kappa ) r^2+(2\mu+\kappa) (\mu+\kappa) (2 \theta \delta) r^{2 \theta+2} }{2\kappa r^2} \notag \\
	&+\frac{\delta \theta (4\mu+2\kappa)r^{2 \theta}}{2\kappa^{\frac{1}{2}} r}+ \delta \theta  |2 \theta-1|(\mu+\kappa)^\frac{1}{2} r^{2 \theta-1}+ 3 \delta^2 \theta^2(\mu+\kappa )^{\frac{1}{2}}  r^{4 \theta-1} \notag \\
	\leq	&C \left( \frac{1}{r}+(1+r^{2\theta})+ r^{2 \theta-1}+ r^{4 \theta-1}\right)\notag \\
	\leq &	C (1+\frac{1}{r}),
\end{align}
which completes the proof.
\hfill
$\Box$


Finally, let us estimate the main term $T(t)$ in \eqref{d71}. Due to 
\[
2 \sin^{2}(tf(r))=1-\cos(2tf(r)),
\]
we have for $t\gg1$, 
\begin{align} \label{d46}
	T(t) &\geq \frac{\omega_{2}}{2}\int_{1/t}^{\varepsilon_0}e^{-r^{2}}\frac{2 \sin^{2}(tf(r))}{f^{2}(r)} rdr\notag\\
	&= \frac{\omega_{2}}{2}\int_{1/t}^{\varepsilon_0}e^{-r^{2}}\frac{r}{f^{2}(r)}dr - \frac{\omega_{2}}{2}\int_{1/t}^{\varepsilon_0}e^{-r^{2}}\frac{r}{f^2(r)}\cos(2tf(r))dr\notag\\
	&=: \frac{\omega_{2}}{2}T_{1}(t) - \frac{\omega_{2}}{2}T_{2}(t).
\end{align}
$T_1(t)$ is estimated by
\begin{align}\label{d47}
	T_{1}(t)&=\int_{\frac{1}{t}}^{\varepsilon_0}e^{-r^{2}}\frac{1+\delta r^{2\theta}}{\mu r^{3}+\kappa r}dr \geq \int_{\frac{1}{t}}^{\varepsilon_0} \frac{e^{-r^{2}}}{(\mu+\kappa)r}dr\geq 
\frac{e^{-\varepsilon_0^2}}{\mu+\kappa} \int_{\frac{1}{t}}^{\varepsilon_0} r^{-1} dr\geq C \log t
\end{align}
\noindent


For $T_{2}(t)$, we develop the integration by parts that idea is inspired from \cite[Proposition A.1]{CT}. Since
\[\cos(2tf(r)) = \frac{1}{2f'(r)t}  \left( \frac{d}{dr}\sin(2tf(r))\right),\]
$T_{2}(t)$ can be written as 
\begin{align}\label{d42}
	T_{2}(t) &= \int_{1/t}^{\varepsilon_0}e^{-r^{2}}\cos(2tf(r))\frac{1+\delta r^{2\theta}}{\mu r^3+\kappa r}dr  =
	\frac{1}{2t}\int_{1/t}^{\varepsilon_0}\frac{e^{-r^{2}}(1+\delta r^{2\theta})}{ f'(r)(\mu r^3+\kappa r)}  ~ d\big(\sin(2tf(r))\big) 
	=: \frac{1}{2t}K(t).
\end{align}
It is obvious that $K(t)$ satisfies
\begin{align} \label{d43}
	K(t) &= \left[\frac{e^{-r^{2}}(1+\delta r^{2\theta})}{ f'(r)(\mu r^3+\kappa r)}  \sin(2tf(r))\right]_{1/t}^{\varepsilon_0} -\int_{1/t}^{\varepsilon_0}\frac{d}{dr}\left(\frac{e^{-r^{2}}(1+\delta r^{2\theta})}{ f'(r)(\mu r^3+\kappa r)} \right)\sin(2tf(r))dr \notag\\
	&=: K_{1}(t) + K_{2}(t).
\end{align}
By Proposition \ref{d40}, it holds that for some universal constants $C > 0$ and $C' > 0$
\begin{equation}\label{d44}
	\vert K_{1}(t)\vert \leq C'\frac{(1+\delta)}{\kappa\varepsilon_{0}} + C'\frac{(1+\delta)}{\kappa (1/t)} \leq C(1+t).
\end{equation}
For the estimate of $ K_{2}(t)$, it follows from  Proposition \ref{d40} that 
\begin{align}
	\vert K_{2}(t)\vert =& \bigg\vert\int_{1/t}^{\varepsilon_0}\bigg( -\frac{e^{-r^{2}}(2r)(1+\delta r^{2\theta})}{ f'(r)(\mu r^3+\kappa r)}+\frac{e^{-r^{2}}(2\delta \theta r^{2\theta-1})}{ f'(r)(\mu r^3+\kappa r)}-\frac{e^{-r^{2}}(1+\delta r^{2\theta})f^{''}(r)}{ \big(f'(r)\big)^2 (\mu r^3+\kappa r)}\notag \\
	&-\frac{e^{-r^{2}}(1+\delta r^{2\theta})(3\mu r^3+\kappa)}{ f'(r)(\mu r^3+\kappa r)^{2}}\bigg)\sin(2tf(r))dr\bigg\vert \notag \\
	\leq& \int_{\frac{1}{t}}^{\varepsilon_0}\frac{2(1+\delta )e^{-r^{2}}}{\kappa |f'(r)| } +\frac{2\delta \theta e^{-r^{2}} r^{2\theta-1}}{\kappa| f'(r)| r} +\frac{(1+\delta)e^{-r^{2}}|f^{''}(r)|}{ \big(f'(r)\big)^2 (\kappa r)}+\frac{(1+\delta)(3\mu +\kappa)e^{-r^{2}}}{ \kappa^{2} |f'(r)| r^{2}} dr\notag\\
	\leq &C \int_{\frac{1}{t}}^{\varepsilon_0} \left(1+r^{2\theta-2} +r^{-1}+r^{-2} \right)dr\notag\\
	 \leq& Ct \label{d45}
\end{align}
with $t\gg 1$.

Combing \eqref{d42}-\eqref{d45},  the estimate for $T_{2}(t)$ is presented by
\begin{equation}\label{d48}
\vert T_{2}(t)\vert \leq C, \quad (t \gg 1). 
\end{equation}
Finally, by \eqref{d71}, \eqref{i-1}, \eqref{d46}, \eqref{d47} and \eqref{d48}, we can obtain the following lemma.\\
\begin{lem}\label{lem3-2} Let $n = 2$, $\kappa>0$,  $\mu \geq0$, $\theta \in (0,2]$ and $\gamma \in (0,1]$. Assume $u_{0} \in L^{2}({\bf R}^{2})$ and $u_{1} \in L^{1,\gamma}({\bf R}^{2})$. Then, it holds that
\[\Vert w(t,\cdot)\Vert^{2} \geq CP^{2}\log t, \quad( t \gg 1).\]
\end{lem}


\section{Estimates from above of the solution: $n = 1,2$}

In this section, by $L^{2}$-regularity and $L^{1}$-regularity of the initial data, the upper bound estimate of $\Vert u(t,\cdot)\Vert$ is given as $t \to \infty$.  

First, we consider the one dimensional case. 
By density arguments, we may assume $u_{0}, u_{1} \in C_{0}^{\infty}({\bf R})$.
It follows from  \eqref{d36} that 
\begin{align} \label{d49}
	\frac{1}{2}\Vert w(t,\cdot)\Vert^{2} &\leq \int_{{\bf R}_{\xi}}\left| \frac{\sin(tf(\vert\xi\vert))}{f(\vert\xi\vert)}\right| ^{2}\vert w_{1}(\xi)\vert^{2}d\xi+
	 \int_{{\bf R}_{\xi}}\vert \cos(tf(\vert\xi\vert))\vert^{2}\vert w_{0}(\xi)\vert^{2}d\xi
		\notag\\
	&\leq \int_{L_{0}}\frac{\sin^{2}(tf(\vert\xi\vert))}{f^{2}(\vert\xi\vert)}\vert w_{1}(\xi)\vert^{2}d\xi + \int_{{\bf R}_{\xi}\setminus L_{0}}\frac{\sin^{2}(tf(\vert\xi\vert))}{f^{2}(\vert\xi\vert)}\vert w_{1}(\xi)\vert^{2}d\xi+ \|w_0(\xi)\|^2 \notag \\
	&=: L_{1}(t) + L_{2}(t)+ \|u_0\|^2, 
\end{align}
where  $L_{0}$ is defined in \eqref{d30} with $n=1$. By \eqref{d50}, $L_{1}(t)$ is estimated by
\begin{equation}\label{15}
L_{1}(t) \leq L^{2}t^{2}\omega_{1}\Vert u_{1}\Vert_{1}^{2}\int_{0}^{\beta(t)}dr = L^{2}t^{2}\omega_{1}\Vert u_{1}\Vert_{1}^{2}\beta(t) \leq C\Vert u_{1}\Vert_{1}^{2}~t
\end{equation}
for $t \gg 1$.

On the other hand, $L_{2}(t)$ can be estimated as follows.
\begin{align}
L_{2}(t) &\leq \int_{{\bf R}_{\xi}\setminus L_{0}}\frac{1+\delta \vert\xi\vert^{2\theta}}{\mu\vert\xi\vert^{4}+\kappa|\xi|^2}\vert w_{1}(\xi)\vert^{2}d\xi\notag\\
&\leq  \int_{\beta(t) \leq \vert\xi\vert \leq \gamma(t)}\frac{1+\delta \vert\xi\vert^{2\theta}}{\mu\vert\xi\vert^{4}+\kappa|\xi|^2} \vert w_{1}(\xi)\vert^{2} d\xi + \int_{\gamma(t) \leq \vert\xi\vert}\frac{1+\delta \vert\xi\vert^{2\theta}}{\mu\vert\xi\vert^{4}+\kappa|\xi|^2} \vert w_{1}(\xi)\vert^{2} d\xi \notag\\
&\leq   \int_{\beta(t) \leq \vert\xi\vert \leq \gamma(t)}\frac{1+\delta \vert\xi\vert^{2\theta}}{\kappa|\xi|^2} \vert w_{1}(\xi)\vert^{2}d\xi + \int_{\gamma(t) \leq \vert\xi\vert}\frac{1+\delta \vert\xi\vert^{2\theta}}{\mu\vert\xi\vert^{4}} \vert w_{1}(\xi)\vert^{2} d\xi \notag\\
&=:  L_{2,1}(t) +  L_{2,2}(t),\quad (t \gg 1)
\end{align}
where $\beta(t)$ is defined in \eqref{d30}  and 
\begin{align} \label{d56}
\gamma(t) := \frac{\delta_{0}}{(\mu+\kappa)^{\frac{1}{2}}  \log t}.
\end{align}
If $\beta(t)\leq |\xi|\leq \gamma(t)$, we can realize $|\xi| \leq 1$ by choosing $t\gg1$. Therefore, we can proceed the estimate of $ L_{2,1}(t)$ as follows:
\begin{equation}\label{key}
	L_{2,1}(t)\leq 2\frac{1+\delta}{\kappa} \|u_1\|_1^2 \int_{\beta(t)}^{\gamma(t)} r^{-2}dr= 2\frac{1+\delta}{\kappa}\|u_1\|_1^2 \left[ -\frac{1}{r}\right]_{\beta(t)}^{\gamma(t) }\leq C\|u_1\|_1^2 ~t.
\end{equation}
For $ L_{2,2}(t)$, we conclude by  $\theta\in [0,2]$ that 
\begin{align} \label{d51}
	 L_{2,2}(t)=& \frac{1}{\mu} \int_{\gamma(t)\leq \vert\xi\vert} |\xi|^{-4}  \vert w_{1}(\xi)\vert^{2}d\xi +  \frac{\delta}{\mu} \int_{\gamma(t)\leq \vert\xi\vert}  |\xi|^{2\theta-4}  \vert w_{1}(\xi)\vert^{2}d\xi\notag\\
	 \leq &  \frac{1}{\mu} \gamma(t)^{-4} \int_{\gamma(t) \leq \vert\xi\vert} \vert w_{1}(\xi)\vert^{2}d\xi+  \frac{\delta}{\mu} ~|\gamma(t)|^{2\theta-4}  \int_{\gamma(t) \leq \vert\xi\vert}\vert w_{1}(\xi)\vert^{2}d\xi\notag \\
	 \leq & C \big((\log t)^4+ (\log t)^{4-2\theta}\big) \|u_1\|^2.
	 \end{align}

The following lemma is a direct consequence of  \eqref{d49}-\eqref{d51}.
\begin{lem}\label{lem3}Let $n = 1$,  $2\geq \theta  > 0$, $\delta>0$, $\mu > 0$ and $\kappa > 0$. Under the assumption $u_0\in  L^{2}({\bf R})$ and $u_{1} \in L^{1}({\bf R})\cap L^{2}({\bf R})$, it holds that
\[\Vert w(t,\cdot)\Vert^{2} \leq C(\Vert u_{1}\Vert^{2} + \Vert u_{1}\Vert_{1}^{2})t + \Vert u_{0}\Vert^{2}, \quad( t \gg 1).\]
\end{lem}


It remains to discuss the upper bound estimate for $n = 2$. We decompose the solution into three parts  in order to get the desired growth estimate.
\begin{align}\label{d57}
\frac{1}{2}\Vert w(t,\cdot)\Vert^{2} &\leq \int_{{\bf R}_{\xi}^{2}}\left| \frac{\sin(tf(\vert\xi\vert))}{f(\vert\xi\vert)}\right| ^{2}\vert w_{1}(\xi)\vert^{2}d\xi +
\int_{{\bf R}_{\xi}^2}\vert \cos(tf(\vert\xi\vert))\vert^{2}\vert w_{0}(\xi)\vert^{2}d\xi \notag\\
&\leq \int_{|\xi|\leq \beta(t)}\frac{\sin^{2}(tf(\vert\xi\vert))}{f^{2}(\vert\xi\vert)}\vert w_{1}(\xi)\vert^{2}d\xi+ \int_{\beta(t)\leq|\xi|\leq 1}\frac{\sin^{2}(tf(\vert\xi\vert))}{f^{2}(\vert\xi\vert)}\vert w_{1}(\xi)\vert^{2}d\xi \notag \\&~~~~+ \int_{|\xi|\geq 1}\frac{\sin^{2}(tf(\vert\xi\vert))}{f^{2}(\vert\xi\vert)}\vert w_{1}(\xi)\vert^{2}d\xi + \|w_0\|^2\notag\\
&=: G_{1}(t) + G_{2}(t)+ G_{3}(t) +\|u_0\|^2, ~~~~~~~~~~(t\gg1),
\end{align}
In the derivation of upper bound estimates, it is enough to assume $u_{0}, u_{1} \in C_{0}^{\infty}({\bf R}^{2})$.\\

We first give the  estimate of $G_{1}(t)$ by \eqref{d50} that 
\begin{equation}\label{}
G_{1}(t) \leq L^{2}t^{2}\omega_{2}\Vert u_{1}\Vert_{1}^{2}\int_{0}^{\beta(t)}rdr = 
L^{2}t^{2}\omega_{2} \Vert u_{1}\Vert_{1}^{2} \left[  \frac{r^2}{2}\right]_0^{\beta(t) }
 \leq C\Vert u_{1}\Vert_{1}^{2}, \quad (t\gg1).
\end{equation}
For the estimate of $G_{2}(t)$, it follows from the fact $|\xi|\leq 1$ and \eqref{d55} that 
\begin{align}
G_{2}(t) &\leq\omega_{2} \Vert  u_{1}\Vert_{1}^{2}\int_{\beta(t)}^{1}\frac{1+\delta r^{2\theta}}{\mu r^3+\kappa r}dr\notag\\ &
 \leq   \omega_{2} \Vert u_{1}\Vert_{1}^{2}\int_{\beta(t)}^{1}\frac{1+\delta }{\kappa r}dr\notag\\ &
 \leq C \Vert u_{1}\Vert_{1}^{2} \left[\log r\right]_{\beta(t)}^{1}\notag\\ &
 \leq C \Vert u_{1}\Vert_{1}^{2} \log t
, \quad (t\gg1).
\end{align}
When $|\xi|\geq 1$, $G_{3}(t)$ can be estimated as 
satisfies
\begin{align}\label{d58}
	G_{3}(t) &\leq \int_{|\xi|\geq 1}\frac{1+\delta r^{2\theta}}{\mu r^4+\kappa r^2} |w_1(\xi)|^2 d\xi\notag\\ &
 \leq \int_{|\xi|\geq 1} \frac{(1+\delta)  r^{2\theta}}{\mu r^4}  |w_1(\xi)|^2 d\xi\notag\\ &
\leq \frac{1+\delta}{\mu} \int_{|\xi|\geq 1}  |w_1(\xi)|^2 d\xi\notag\\ &
	\leq C \|u_1\|^2
	, \quad (t\gg1),
\end{align}
where we used the fact $\theta \in (0,2]$.\\

By the above arguments from \eqref{d57} to \eqref{d58}, we conclude the following lemma.\\

\begin{lem}\label{lem44}Let $n = 2$, $\delta  > 0$, $2\geq \theta  > 0 $, $\mu > 0$ and $\kappa > 0$. Under the assumption $u_{0} \in L^{2}({\bf R}^{2})$ and $u_{1} \in L^{1}({\bf R}^{2})\cap L^{2}({\bf R}^{2})$, it holds that
\[\Vert w(t,\cdot)\Vert^{2} \leq C(\Vert u_{1}\Vert^{2} +\Vert u_{1}\Vert_{1}^{2})\log t + \Vert u_{0}\Vert^{2},\quad(t \gg 1).\]
\end{lem}

\par
\vspace{0.5cm}
Finally, by the Plancherel theorem, we can see that Theorems \ref{theorem1} and \ref{theorem2} can be demonstrated by  Lemmas \ref{lem1}, \ref{lem3-2}, \ref{lem3} and \ref{lem44}.\\


\section{$L^{2}$-upper bound for the case $n \geq 3$}

In this section, we mention the boundedness of the $L^2$ norm of the solution itself in the case of space $3$ dimensions or more.\\

At first, we can proceed the estimate by \eqref{d36} that
\begin{align}
	\frac{1}{2}\Vert w(t,\cdot)\Vert^{2} &\leq \int_{{\bf R}_{\xi}^{n}}\left| \frac{\sin(tf(\vert\xi\vert))}{f(\vert\xi\vert)}\right|  ^{2}\vert w_{1}(\xi)\vert^{2}d\xi  +\int_{{\bf R}_{\xi}^{n}}\vert\cos(tf(\vert\xi\vert))\vert ^{2}\vert w_{0}(\xi)\vert^{2}d\xi  \notag\\
	&\leq \int_{\vert\xi\vert \leq 1}\frac{\sin^{2}(tf(\vert\xi\vert))}{f^{2}(\vert\xi\vert)}\vert w_{1}(\xi)\vert^2 d\xi + \int_{\vert\xi\vert \geq 1}\frac{\sin^{2}(tf(\vert\xi\vert))}{f^{2}(\vert\xi\vert)}\vert w_{1}(\xi)\vert^{2}d\xi + \int_{{\bf R}_{\xi}^{n}}\vert w_{0}(\xi)\vert^{2}d\xi   \notag\\
	&=: M_{1}(t) + M_{2}(t)+\|w_0\|^2. \label{d72}
\end{align}
Because of $n\geq 3$ and $\theta \in[0,2]$, $M_1(t)$ and $M_2(t)$ can be estimated by
\begin{align}\label{d73}
	M_{1}(t) &\leq \int_{|\xi|\leq 1}\left| \frac{\sin(tf(|\xi|))}{f(\vert\xi\vert)}\right|^2\vert w_{1}(\xi)\vert^{2}d\xi\notag\\&
	 \leq \omega_{n}\Vert u_{1}\Vert_{1}^{2}\int_{0}^{1}\frac{(1+\delta r^{2\theta})r^{n-1}}{\mu r^4 +\kappa r^2}dr \notag\\
	& \leq \omega_{n} \frac{1+ \delta}{\kappa} \Vert u_{1} \Vert_{1}^{2} \int_{0}^{1} r^{n-3}dr\notag\\&\leq C\Vert u_{1}\Vert_{1}^{2},
\end{align}

\begin{align}\label{d74}
	M_{2}(t) &\leq \int_{|\xi| \geq 1 }\frac{1+\delta |\xi|^{2\theta}}{\mu |\xi|^4 +\kappa|\xi|^2}|\omega_{1}(\xi)|^{2} d\xi \notag\\
	&\leq 
	\int_{|\xi| \geq 1 }\frac{(1+\delta)}{\mu |\xi|^{4-2\theta} }|\omega_{1}(\xi)|^{2} d\xi\notag\\&
	 \leq 
	\frac{\delta+1}{\mu} \int_{|\xi| \geq 1} |\omega_{1}(\xi)|^{2}  d\xi\notag\\&
	\leq C\Vert u_{1}\Vert^{2}.
\end{align}

Thus, by \eqref{d72}, \eqref{d73} and \eqref{d74}, we can obtain the following desired estimate.
\begin{pro}\label{lem4}Let $n \geq 3$, $0 \leq \theta \leq 2$, $\mu > 0$, $\kappa >  0$, $u_{0} \in L^{2}({\bf R}^{n})$ and $u_{1} \in L^{1}({\bf R}^{n})\cap L^{2}({\bf R}^{n})$. Then, it holds that
\[\Vert u(t,\cdot)\Vert^{2} = \Vert w(t,\cdot)\Vert^{2} \leq C(\|u_0\|^2+\Vert u_{1}\Vert^{2} + \Vert u_{1}\Vert_{1}^{2}), \quad ( t \gg 1).\]
\end{pro}
\begin{rem}{\rm From the proof of \eqref{d74} it can be seen that the $L^{2}$-bound estimate of the solution itself in the high-frequency band is essentially brought about by the effect of the plate equation with $\mu > 0$.} 

\end{rem}

\section{Remark on Hardy-type inequalities in low dimensions}

In this section, how Hardy's inequality does not hold in the case of $1$, $2$ dimensional whole spaces, and we give a note on this. As is well-known, the so-called Hardy's or Pitt's inequality has a form: 
\begin{equation}\label{1}
\Vert\frac{u}{d(\cdot)}\Vert_{L^{2}({\bf R}^{n})} \leq C^{*}\Vert\nabla u\Vert_{L^{2}({\bf R}^{n})}\quad (\forall u \in H^{1}({\bf R}^{n}))
\end{equation}
with some universal constant $C^{*} > 0$. Here $n \geq 3$ (see \cite{Be}), and 
\begin{equation}\label{2}
d(x) := \vert x\vert \quad \textstyle{or} \quad 1+\vert x\vert.
\end{equation}
In the case when ${\bf R}^{n}$ is replaced by a smooth exterior domain $\Omega \subset {\bf R}^{n}$ ($0 \notin \bar{\Omega}$) with compact boundary, the inequality \eqref{1} still holds true even for $n = 2$ (see \cite{DS}):
\begin{equation}\label{3}
\Vert\frac{u}{d(\cdot)}\Vert_{L^{2}(\Omega)} \leq C^{*}\Vert\nabla u\Vert_{L^{2}(\Omega)}\quad (\forall u \in H_{0}^{1}(\Omega))
\end{equation}
with some universal constant $C^{*} > 0$, where
\[d(x) := \left\{
   \begin{array}{ll}
     \vert x\vert& 
       \qquad n \geq 3, \\[0.2cm]
      \log(B\vert x\vert) & \qquad n = 2,
    \end{array} \right.\]
and $B > 0$ is a constant satisfying 
\[B\vert x\vert \geq 2\quad (\forall x \in \Omega).\]
   
The question then naturally arises whether inequalities of the form \eqref{3} hold in the case of a two-dimensional ``whole space", or whether they hold in the case of a one-dimensional whole space. In this paper, we use the properties of solutions to the Rosenau equation just obtained in Theorems \ref{theorem1} and \ref{theorem2} to provide an answer to this question. \\

For simplicity, we restrict our discussion to the following case of a weight function $w(x)$ in two spatial dimensions. The considerations obtained below are not new, having already been discussed in \cite{MOW}, but they are meant as a prototype to the weight function $w(x)$ in space $1$-D or another type of weight function. The Hardy-type inequality does not hold in $1$ or $2$ dimensions. It is an interesting point of view that the non-sufficiency of Hardy-type inequalities independent of the equations can be confirmed by the nature of the solution of a particular equation (here Rosenau equation).
For that purpose, we fix $n = 2$ and restrict ourselves to the case of a weight function $w(x)$ (as one example of a Hardy-type inequality that does not hold) defined by\\
{\bf (A-1)}\,\,\,$w(x) := (1+\log(1+\vert x\vert))(1+\vert x\vert)$.\\

Our statement reads as follows. 
\begin{pro}\label{theorem3}
Let $n = 2$ and assume {\bf (A-1)}. Then, the following inequality never holds.
\begin{equation}\label{4}
\Vert\frac{u}{w(\cdot)}\Vert_{L^{2}({\bf R}^{2})} \leq C^{*}\Vert\nabla u\Vert_{L^{2}({\bf R}^{2})}\qquad (\forall u \in H^{1}({\bf R}^{2})).
\end{equation}
\end{pro}
\begin{rem}{\rm One can check Proposition \ref{theorem3} by replacing $w(x)$ by $w(x) := (1+\vert\log\vert x\vert\vert)(1+\vert x\vert)$, much less the following, which is another example where Hardy's inequality does not hold: $w(x) := \vert x\vert(1+\vert\log\vert x\vert\vert)$ (for these examples, see \cite{MOW}). Additionally, if we take $w(x) = \vert x\vert$ in the $1$-dimensional case, it still does not hold. Of course, other cases can be checked in the following ways developed in subsection 5.1.}
\end{rem}
\begin{rem}{\rm By taking $w(x) = 1$ in \eqref{4} formally, we also see that the Poincar\'e-type inequality does not hold in the case of $n = 1,2$.}
\end{rem}

\subsection{Proof of Proposition \ref{theorem3}}

In this subsection, we prove Proposition \ref{theorem3} by using several properties of the solution itself to the Cauchy problem \eqref{d1}-\eqref{d2} of the Rosenau equation with $\theta = 1$ and $n = 2$:
\begin{numcases}
  ~u_{tt} -\delta \Delta u_{tt} +\mu  \Delta^2u-\kappa  \Delta u = 0,\ \ \ (t,x)\in (0,\infty)\times {\bf R}^{2},\label{f1}\\
 u(0,x) = 0, \quad  u_{t}(0,x)= u_{1}(x),\ x\in{\bf R}^{2}, \label{f2}
\end{numcases}
where $\delta >0$,~$\mu>0$, $\kappa>0$. Here, we choose initial velocity $u_{1} \in C_{0}^{\infty}({\bf R}^{2})$, and 
$$\int_{{\bf R}^{2}}u_{1}(x) dx \ne 0.$$
Note that we just have chosen $u_{0}(x) = 0$. Under this assumption, the standard regularity argument shows that the solution $u(t,x)$ to problem \eqref{f1}-\eqref{f2} is sufficiently smooth in the time-space direction (cf. \cite{ikawa}). Since $u_{1} \in L^{1,1}({\bf R}^{2})$, by applying Theorem \ref{theorem2}, for the corresponding (smooth) solution $u(t,x)$ of the problem \eqref{f1} and \eqref{f2} we see
\begin{equation}\label{f3}
\lim_{t \to \infty}\Vert u(t,\cdot)\Vert = +\infty.
\end{equation}

Next, suppose that the inequality \eqref{4} is correct with a constant $C^{*} > 0$ and a function $w(x)$ satisfying {\bf (A-1)}. While, for the smooth solution $u(t,x)$ to problem \eqref{f1}-\eqref{f2} we set
\[v(t,x) := \int_{0}^{t}u(s,x)ds.\]
Then, the function $v(t,x)$ satisfies
\begin{numcases}
  ~v_{tt} -\delta \Delta v_{tt} +\mu  \Delta^2v-\kappa  \Delta v = u_{1},\ \ \ (t,x)\in (0,\infty)\times {\bf R}^{2},\label{f4}\\
 v(0,x) = 0, \quad  v_{t}(0,x) = 0,\ x\in{\bf R}^{2}. \label{f5}
\end{numcases}

This idea comes from \cite{IM} again. By multiplying both sides of \eqref{f4} by $v_{t}(t,x)$, and integrating over $[0,t]\times {\bf R}^{2}$ it follows that
\begin{equation}\label{8}
\frac{1}{2}\Vert v_{t}(t,\cdot)\Vert^{2} + \frac{\delta}{2}\Vert\nabla v_{t}(t,\cdot)\Vert^{2} + \frac{\mu}{2}\Vert\Delta v(t,\cdot)\Vert^{2} + \frac{\kappa}{2}\Vert\nabla v(t,\cdot)\Vert^{2} = \int_{{\bf R}^{2}}u_{1}(x)v(t,x)dx.
\end{equation}
Now, let us estimate the last term of \eqref{8} by using the inequality \eqref{4} whose formation is assumed to be.
Indeed, for any $\varepsilon > 0$ it follows that
\begin{align}
	\left| \int_{{\bf R}^{2}}u_{1}(x)v(t,x)dx\right|  &\leq \int_{{\bf R}^{2}}w(x)\vert u_{1}(x)\vert\frac{\vert v(t,x)\vert}{w(x)}dx \notag\\
	&\leq C_{\varepsilon}\int_{{\bf R}^{2}}w(x)^{2}\vert u_{1}(x)\vert^{2}dx + \varepsilon\int_{{\bf R}^{2}}\frac{\vert v(t,x)\vert^{2}}{w(x)^{2}}dx \notag\\
	&\leq C_{\varepsilon}\int_{{\bf R}^{2}}w(x)^{2}\vert u_{1}(x)\vert^{2}dx + \varepsilon(C^{*})^{2}\Vert\nabla v(t,\cdot)\Vert^{2}
\end{align} \label{9}
with some constant $C_{\varepsilon} > 0$ depending on $\varepsilon > 0$. Thus, from \eqref{8} and \eqref{9} we find that
\begin{equation}\label{10}
\frac{1}{2}\Vert v_{t}(t,\cdot)\Vert^{2} + \left(\frac{\kappa}{2}-\varepsilon(C^{*})^{2}\right)\Vert\nabla v(t,\cdot)\Vert^{2} \leq C_{\varepsilon}\int_{{\bf R}^{2}}w(x)^{2}\vert u_{1}(x)\vert^{2}dx < +\infty.
\end{equation}
By choosing $\varepsilon > 0$ in \eqref{10} sufficiently small, one can arrive at the following crucial estimate: 
\begin{equation}\label{11}
\frac{1}{2}\Vert v_{t}(t,\cdot)\Vert^{2} \leq C_{\varepsilon}\int_{{\bf R}^{2}}w(x)^{2}\vert u_{1}(x)\vert^{2}dx.
\end{equation}
Since $v_{t} = u$, one has the bounded estimate of the solution itself:
\begin{equation}\label{12}
\Vert u(t,\cdot)\Vert^{2} \leq 2C_{\varepsilon}\int_{{\bf R}^{2}}w(x)^{2}\vert u_{1}(x)\vert^{2}dx < +\infty.
\end{equation}
\noindent
The estimate \eqref{12} contradicts \eqref{f3} in the case of $n = 2$. This shows the validity of the statement of Proposition \ref{theorem3}.\\
Incidentally, the above proof can be expressed in the form of a motto: "If the inequality \eqref{4} is correct, then the $L^{2}$-boundedness of the solution itself can be derived".
\hfill
$\Box$
\begin{rem}{\rm When we choose $w(x) = \vert x\vert(1+\vert\log\vert x\vert\vert)$ as the weight function, even if $u_{1} \in C_{0}^{\infty}({\bf R}^{2})$, it is not trivial to check $\int_{{\bf R}^{2}}w(x)^{2}\vert u_{1}(x)\vert^{2}dx < +\infty$. In this case it suffices to start with the proof by choosing $u_{1} \in C_{0}^{\infty}({\bf R}^{2})$ satisfying $u_{1}(x) = 0$ near $x = 0$ and $\int_{{\bf R}^{2}}u_{1}(x)dx \ne 0$. This is feasible.}
\end{rem}
\begin{rem}{\rm 
By basing on the similar argument above and using the growth results derived in \cite{FE-ike} about the plate equation, one can also discuss the failure of the so-called Hardy-Rellich inequality (cf. \cite{DH}) in the case of $n = 1,2,3,4$.
\begin{cor}\,Let $n \geq 5$. Then there exists a constant $C^{*} > 0$ such that
\[\left\|  \frac{u}{\vert x\vert^{2}}\right\|  \leq C^{*}\Vert\Delta u\Vert\]
for all $u \in H^{2}({\bf R}^{n})$.
\end{cor}
}
\end{rem}
\vspace{0.5cm}
\noindent{\em Acknowledgement.}
\smallskip
This paper was written during Xiaoyan Li's stay as an overseas researcher at Hiroshima University from 12 December, 2022 to 11 December, 2023 under Ikehata's supervision as a host researcher. The work of the first author (Xiaoyan Li) was financially supported in part by Chinese Scholarship Council (Grant No. 202206160071). 
The work of the second author (Ryo Ikehata) was supported in part by Grant-in-Aid for Scientific Research (C) 20K03682 of JSPS. 
\vspace{0.5cm}

\noindent{\bf Declarations}
\smallskip

\noindent{\bf Data availability}
\smallskip
Data sharing not applicable to this article as no datasets were generated or analysed during the current study.

\noindent{\bf Conflict of interest}
\smallskip
The authors declare that they have no conflict of interest.

\section{Appendix}

{\small The well-posedness as in Proposition \ref{prp1} of the solution for \eqref{d1}-\eqref{d2}, which is more or less known from previous studies by \cite{H2016} and \cite{mh2019}, is included here as an appendix for the reader's convenience. In particular, we follow the idea due to \cite{H2016} and \cite[Theorem A.1]{mh2019}. In fact, Michihisa \cite{mh2019} treated the strongly damped Rosenau equation: for $\nu > 0$
\[u_{tt}-\Delta u -\nu\Delta u_{t}+\Delta^{2}u + \Delta^{2}u_{tt} = 0.\]
The following proof is a slight modification of \cite[Theorem A.1]{mh2019}.\\

Considering the energy space $\mathcal{H}:= H^2({\bf R}^n)\times  H^\theta({\bf R}^n)$ ($0 \leq \theta \leq 2$), the inner product in $\mathcal{H}$ is defined by 
\[
\big([u^1,v^1], ~[u^2,v^2] \big)_\mathcal{H}=(u^1, u^2)_{H^2}+(v^1, v^2)_{H^\theta}.
\]
Setting $v=u_t$ and $A=-\Delta$, it follows from \eqref{d1} that 
$$(I+\delta A^\theta )v_t=-(\mu A^2+\kappa A+I )u+ u,$$
where $I$ is the identity operator.
Then the problem \eqref{d1}-\eqref{d2} can be written as 
\begin{equation} \label{d18}
	\left\{
	\begin{aligned}
		&\frac{d}{dt}U=\mathcal{A}U+ \mathcal{B}U,\\
		&U(0,x) =U_0,
	\end{aligned}
	\right.
\end{equation}
where 
\begin{align*}
\mathcal{A} =	\begin{pmatrix}
	0& I \\ -P& 0
\end{pmatrix},~~~
P = (I+\delta A^{\theta})^{-1}(\mu A^2+\kappa A +I), ~~~~~~~\notag \\
\mathcal{B} =	\begin{pmatrix}
	0& 0 \\ (I+\delta A^\theta )^{-1}& 0
\end{pmatrix},~~~
U=	\begin{pmatrix}
	u\\v
\end{pmatrix}\in \mathcal{H},~~~
U_0=	\begin{pmatrix}
	u_0\\u_1
\end{pmatrix} \in \mathcal{H}.
\end{align*}
Here $\mathcal{D}(P)$ is defined by 
\begin{align}
	\mathcal{D}(P)
	=\biggl\{ &u\in H^2({\bf R}^n): \text{There exists} ~y=y_u\in H^\theta({\bf R}^n) 
	~\text{such that} ~ \notag \\  
	&(u,\phi)+\kappa(A^{\frac{1}{2}}u, A^{\frac{1}{2}}\phi)+\mu(Au, A\phi) =(y,\phi)+\delta(A^{\frac{
		\theta}{2}}y, A^{\frac{\theta}{2}}\phi ),~~\forall \phi \in  H^2({\bf R}^n) \biggr\}
 . \label{d4}
\end{align}

Note  that $\mathcal{D}(P)$ is not empty because $0\in \mathcal{D}(P)$ when  $y=0$.  If $u\in \mathcal{D}(P)$, then $y\in H^\theta({\bf R}^n)$ such that \eqref{d4} holds is uniquely determined. Otherwise, at least there exists  $y^1$ and $y^2 \in H^\theta({\bf R}^n)$ satisfying 
\begin{align}
(u,\phi)+\kappa(A^{\frac{1}{2}}u, A^{\frac{1}{2}}\phi)+\mu(Au, A\phi) =(y_1,\phi)+\delta(A^{\frac{
		\theta}{2}}y_1, A^{\frac{\theta}{2}}\phi )\label{d5} \\ 
(u,\phi)+\kappa(A^{\frac{1}{2}}u, A^{\frac{1}{2}}\phi)+\mu(Au, A\phi) =(y_2,\phi)+\delta(A^{\frac{
		\theta}{2}}y_2, A^{\frac{\theta}{2}}\phi ) \label{d6}
\end{align}
for all $\phi \in  H^2({\bf R}^n)$. Taking $\bar{y}=y^1-y^2\in  H^\theta ({\bf R}^n)$, \eqref{d5} and \eqref{d6}  imply
\begin{equation}\label{d64}
(\bar{y}, \phi)+\delta(A^{\frac{\theta}{2}}\bar{y}, A^{\frac{\theta}{2}}\phi )=0.
\end{equation}
Owing that $C_0^\infty({\bf R}^n)$ is dense in $ H^2({\bf R}^n)$, \eqref{d64} is still correct for all $\phi \in C_0^\infty({\bf R}^n)$. Because $C_0^\infty({\bf R}^n)$ is also dense in $ H^\theta({\bf R}^n)$, we choose a sequence $\{ y_k \} \in C_0^\infty({\bf R}^n)$ such that $\lim\limits_{k\rightarrow \infty }y_k \rightarrow \bar{y}$ in $ H^\theta({\bf R}^n)$, that is
\begin{align} \label{d65}
	\|y_k-\bar{y}\|_{H^\theta}^2=\|\bar{y}\|_{H^\theta}^2-2(\bar{y}, y_k)_{H^\theta}+\|y_k\|_{H^\theta}^2\rightarrow0,\quad (k\rightarrow \infty).
\end{align}
It follows from $ \big| \|y_k\|_{H^\theta}-\|\bar{y}\|_{H^\theta} \big| \leq 	\|y_k-\bar{y}\|_{H^\theta} \rightarrow0$ that 
\begin{align} \label{d66}
	\|y_k\|_{H^\theta} \rightarrow \|\bar{y}\|_{H^\theta} ,\quad (k\rightarrow \infty).
\end{align}
Combining \eqref{d65} and \eqref{d66} yields 
\begin{align}\label{d68}
	\lim\limits_{k\rightarrow \infty }(\bar{y}, y_k)_{H^\theta}=\|\bar{y}\|_{H^\theta}^2.
\end{align}
By \eqref{d64} and the definition of inner product in $H^\theta({\bf R}^n)$, we have 
\begin{align} \label{d67}
0=(\bar{y}, y_k)+\delta(A^{\frac{\theta}{2}}\bar{y}, A^{\frac{\theta}{2}}y_k )=(\bar{y}, y_k)_{H^\theta}
\end{align}
We conclude by \eqref{d68} and \eqref{d67} that 
\begin{align}
		\lim\limits_{k\rightarrow \infty }(\bar{y}, y_k)_{H^\theta}=\|\bar{y}\|_{H^\theta}^2=0,
\end{align}
which implies $\bar{y}=0$ in $H^\theta({\bf R}^n)$, that is, $y^{1} = y^{2}$.

The above arguments show that the linear operator $P: u\rightarrow y=y_u $ is well defined for each $u \in \mathcal{D}(P)$.\\

Next we demonstrate the fact 
\begin{equation}\label{d17}
\mathcal{D}(P)=H^{4-\theta}({\bf R}^n)	
\end{equation}
for $0\leq\theta\leq 2$ by Lemma \ref{d7} and \ref{d8} below.
\begin{lem} \label{d7}
If $0\leq\theta\leq 2$, it holds that $\mathcal{D}(P)\subset H^{4-\theta}({\bf R}^n)$ and there exists a constant $C>0$ such that 
\[
\|u\|_{H^{4-\theta}}\leq C\|Pu\|_{H^\theta}\]
for all $u\in \mathcal{D}(P)$.
\end{lem}
\noindent{ \it{Proof.}}
For any given $u\in \mathcal{D}(P)$, there exists $y=y_u\in H^\theta$ such that
\begin{equation}\label{d9}
 (u,\phi)+\kappa(A^{\frac{1}{2}}u, A^{\frac{1}{2}}\phi)+\mu(Au, A\phi) =(y,\phi)+\delta(A^{\frac{
		\theta}{2}}y, A^{\frac{\theta}{2}}\phi ),~~\forall \phi \in  H^2({\bf R}^n).
\end{equation}
Taking  the Fourier transform on both sides of \eqref{d9}, we have
\[
\int_{{\bf R}^n_\xi} 
(1+\kappa |\xi|^2+\mu|\xi|^4 ) \hat{u} \bar{\hat{\phi}} ~d\xi=\int_{{\bf R}^n_\xi} (1+\delta |\xi|^{2\theta })\hat{y} \bar{\hat{\phi}} ~ d\xi.
\]
Due to the arbitrariness of $\phi$, we can obtain $(1+\kappa |\xi|^2+\mu|\xi|^4 ) \hat{u}=(1+\delta |\xi|^{2\theta })\hat{y}$, that is 

\begin{equation}\label{d10}
	\hat{y}=\widehat{Pu}=\frac{1+\kappa |\xi|^2+\mu|\xi|^4 }{1+\delta |\xi|^{2\theta }} \hat{u}.
\end{equation}
Calculating the $L^2$-norm about $\xi$ on  both sides of \eqref{d10} yields
\begin{equation}\label{d11}
\int_{{\bf R}^n_\xi} (1+\delta |\xi|^{2\theta })|\hat{y} |^2
~d\xi=\int_{{\bf R}^n_\xi}  
\frac{(1+\kappa |\xi|^2+\mu|\xi|^4 )^2  }{(1+\delta |\xi|^{2\theta })} 
|\hat{u}|^2 ~ d\xi.
\end{equation}
Note that $1+|\xi|^{2(4-\theta)}$  can be controlled by $\frac{(1+\kappa |\xi|^2+\mu|\xi|^4 )^2  }{(1+\delta |\xi|^{2\theta })} $.\\
In fact, if we set $r := \vert\xi\vert$ and
\[h(r) := \frac{(1+\kappa |\xi|^2+\mu|\xi|^4 )^2  }{(1+\delta |\xi|^{2\theta })} = \frac{(1+\kappa r^2+\mu r^4 )^2  }{(1+\delta r^{2\theta })},\]
we easily see that
\[\lim_{r \to 0}\frac{h(r)}{1+r^{2(4-\theta)}} = 1,\]
and
\[\lim_{r \to +\infty}\frac{h(r)}{1+r^{2(4-\theta)}} = \frac{\mu^{2}+2\kappa\mu}{1+\delta} := c_{0} > 0.\]
(At this stage, it must be $\mu > 0$.) Thus, there are positive constants $r_{0} \ll 1$ and $r_{\infty} \gg 1$ satisfying $0 < r_{0} < 1 < r_{\infty}$ such that
\[\frac{1}{2} \leq \frac{h(r)}{1+r^{2(4-\theta)}}\quad (0 \leq \forall r \leq r_{0}),\]
\[\frac{c_{0}}{2} \leq \frac{h(r)}{1+r^{2(4-\theta)}}\quad (r_{\infty} \leq \forall r ).\]
Furthermore, because of the continuity of the function $r \mapsto \frac{h(r)}{1+r^{2(4-\theta)}}$ there is a constant $m > 0$ such that
\[m \leq \frac{h(r)}{1+r^{2(4-\theta)}}\quad (r_{0} \leq \forall r \leq r_{\infty}).\]
Set
\[M := \min\{\frac{1}{2}, m, \frac{c_{0}}{2}\}.\]
Then, one can get the lower bound estimate:
\[\frac{(1+\kappa |\xi|^2+\mu|\xi|^4 )^2  }{(1+\delta |\xi|^{2\theta })} \geq M(1+\vert\xi\vert^{2(4-\theta)})\quad (\forall \xi \in {\bf R}_{\xi}^{n}).\]

Therefore, by \eqref{d11}, we have 
\begin{equation}\label{d12}
	\int_{{\bf R}^n_\xi} (1+\delta |\xi|^{2\theta })|\hat{y} |^2
	~d\xi\geq M \int_{{\bf R}^n_\xi}  
(	1+|\xi|^{2(4-\theta)})
	|\hat{u}|^2 ~ d\xi,
\end{equation}
which implies 
\begin{equation}\label{d13}
	\|y\|_{H^\theta}^2\geq M \|u\|_{H^{4-\theta}}^2.
\end{equation}
Thanks to $y=Pu\in H^\theta({\bf R}^n)$, we can obtain $u\in  H^{4-\theta}({\bf R}^n)$.
\hfill
$\Box$

\begin{lem} \label{d8}
If $0\leq\theta\leq 2$, it holds that $ H^{4-\theta}({\bf R}^n) \subset \mathcal{D}(P)$.
\end{lem}
\noindent{ \it{Proof.}}
For any given $u\in H^{4-\theta}({\bf R}^n)$, we define the linear  functional $G_u: H^\theta({\bf R}^n)\rightarrow {\bf R}$ as
\begin{equation}\label{}
	\left\langle G_u, \psi\right\rangle =(u, \psi)+\kappa(Au, \psi)+\mu(A^{2-\frac{\theta}{2}}u, A^\frac{\theta}{2}\psi) 
	,~~~~  \forall \psi \in H^\theta({\bf R}^n).	
\end{equation}
It is not difficult to see that the linear  functional $G_u$ is bounded. In fact, it follows from $u\in H^{4-\theta}({\bf R}^n)$ and $ \psi \in H^\theta({\bf R}^n)$ that
\begin{align*}\label{}
\left| 	\left\langle G_u, \psi \right\rangle \right|  &\leq \left| (u, \psi) \right| +\kappa \left| (Au, \psi)\right| +\mu\left| (A^{2-\frac{\theta}{2}}u, A^\frac{\theta}{2}\psi) \right| \\
&\leq \|\hat{u}\|~\|\hat{\psi}\|+\kappa \|\:|\xi|^{2}\hat{u}\|~\| \hat{\psi}\|+\mu\|\: |\xi|^{4-\theta}\hat{u}\|~ \|\: |\xi|^{\theta}\hat{\psi}\|    \\
&\leq(1+\kappa+\mu)\|u\|_{H^{4-\theta}}\|\psi\|_{H^\theta}.
\end{align*}
Therefore, we have $G_u\in H^{-\theta}{({\bf R}^n)}$.

We also need to define a bilinear functional $a(\cdot, \cdot)$:~$H^\theta({\bf R}^n)\times H^\theta({\bf R}^n)\rightarrow {\bf R} $ as follows
\begin{equation}\label{key}
	a(\psi_1, \psi_2)=(\psi_1, \psi_2)+\delta (A^\frac{\theta}{2}\psi_1, A^\frac{\theta}{2}\psi_2),~~~~~~~~\forall ~ \psi_1, \psi_2 \in  H^\theta({\bf R}^n).
\end{equation}
Note that the bilinear functional $a(\cdot, \cdot)$ are  continuous and coercive because 
\[
|a(\psi_1, \psi_2)|\leq \|\hat{\psi}_1\|~\|\psi_2\|+\|\;|\xi|^\theta\hat{\psi}_1\|~\| \;|\xi|^\theta\hat{\psi}_2\|\leq 2\|\psi_1\|_{H^\theta}\|\psi_2\|_{H^\theta},
\]
and 
\[
a(\psi_1, \psi_1)=\|\psi_1\|_{H^\theta}^2.
\]
Thus, by the Lax-Milgram theorem, for any given $u \in  H^{4-\theta}({\bf R}^n)$, there exists a unique $y = y_{u} \in  H^\theta({\bf R}^n)$ such that \eqref{d14} below holds:
\begin{equation}\label{d14}
\left\langle G_u, \psi\right\rangle = (y, \psi), \forall u \in  H^{4-\theta}({\bf R}^n).
\end{equation}

Next it is necessary to modify the regularity of $\psi$. Due to $0\leq \theta\leq 2$ and $H^2({\bf R}^n)\subset H^\theta({\bf R}^n)$, we can require test functions $\psi \in H^2({\bf R}^n)$ such that 
\begin{equation}\label{d15}
	(Au, \psi)=(A^\frac{1}{2}u, A^\frac{1}{2}\psi),~~~~~~~~~(A^{2-\frac{\theta}{2}}u, A^\frac{\theta}{2}\psi)=(Au, A\psi).
\end{equation}
Therefore, under the condition $\psi \in H^2({\bf R}^n)$, substituting \eqref{d15} to \eqref{d14} yields
\begin{equation}\label{d16}
	(u, \psi)+\kappa(A^\frac{1}{2}u, A^\frac{1}{2}\psi)+\mu(Au, A\psi)=(y, \psi)+\delta(A^\frac{\theta}{2}y, A^\frac{\theta}{2}\psi).
\end{equation}

This implies $u \in \mathcal{D}(P)$, which completes the proof of  the statement $H^{4-\theta}({\bf R}^n) \subset \mathcal{D}(P)$.
\hfill
$\Box$
\vspace{0.5cm}

The fact \eqref{d17} can be guaranteed by Lemmas \ref{d7} and \ref{d8}. Basing on this fact one can prove the crucial lemma.
\begin{lem} \label{d29}
The linear operator 
\begin{gather*}
	\mathcal{A} =	\begin{pmatrix}
		0& I \\ -P& 0
	\end{pmatrix}: 
	\mathcal {D}(	\mathcal {A}) \rightarrow \mathcal{H}
\end{gather*}
generates a strongly continuous contraction semigroup $T(t)$ on $\mathcal {H}$, where 
$$\mathcal{D}(\mathcal{A}) = H^{4-\theta}({\bf R}^n)\times H^2({\bf R}^n) = \mathcal{D}(P)\times H^2({\bf R}^n).$$
\end{lem}
\noindent{ \it{Proof.}}
First, we prove that the operator $\mathcal{A}$ is dissipative. For any given $U=[u,v]^T \in  \mathcal{D}(\mathcal{A})$, it follows from $\eqref{d10}$ that
\begin{align}
	    \big( \mathcal{A}U, U   \big)_\mathcal{H}&=\big( [v, -Pu],~[u,v] \big)_\mathcal{H} \notag \\
	&=(v,u)_{H^2}-(Pu,v)_{H^\theta}\notag\\
	&=\int_{{\bf R}^n}(1+ \kappa |\xi|^2+\mu|\xi|^4 ) \hat{v}\bar{\hat{u}}~d\xi-
	\int_{{\bf R}^n} (1+\delta |\xi|^{2\theta})\widehat{Pu} \bar{\hat{v}}~d\xi  \notag\\
	&=\int_{{\bf R}^n} (1+ \kappa |\xi|^2+\mu|\xi|^4 ) ( \hat{v}\bar{\hat{u}}- \hat{u}\bar{\hat{v}})~
d\xi\notag\\ 
	&=2i \int_{{\bf R}^n}  (1+ \kappa |\xi|^2+\mu|\xi|^4 ) \text{Img}(\hat{u}\bar{\hat{v}})~d\xi,
\end{align}
which implies
\[{\text Re}\big( \mathcal{A}U, U   \big)_\mathcal{H} = 0.\]

Next we prove that $ \mathcal{I}-\mathcal{A}$ is surjective, that is, for any given$[f,g]^T \in \mathcal{H}$, we need to find $[u,v]^T \in \mathcal{D}(\mathcal{A})$ such that 
\begin{gather} \label{d20}
		(\mathcal{I}-\mathcal{A})	\begin{pmatrix}
		u\\v
	\end{pmatrix}
	=	\begin{pmatrix}
		f\\g
	\end{pmatrix},
\end{gather}
which is equivalent to
	\begin{align}
		&u-v=f,   \label{d21}\\
		&Pu+v=g.   \label{d22}
\end{align}
Substituting \eqref{d21} into \eqref{d22} leads to
\begin{equation} \label{d69}
	Pu+u=f+g \in H^{\theta}({\bf R}^{n}).
\end{equation}

Define the bilinear  functional on  $H^2({\bf R}^n)$ as follows
\[
~~~~~~~B(\psi_1, \psi_2)=(\psi_1, \psi_2)_{H^2}+(\psi_1, \psi_2)_{H^\theta},~~~~\forall ~\psi_1, \psi_2 \in H^2({\bf R}^n).
\]
It is not difficult to see that the bilinear functional  $B(\psi_1, \psi_2)$ is continuous and coercive, that is,
\[
|B(\psi_1, \psi_2)|\leq \|\psi_1\|_{H^2}\|\psi_2\|_{H^2}+\|\psi_1\|_{H^\theta}\|\psi_2\|_{H^\theta}\leq C\|\psi_1\|_{H^2}\|\psi_2\|_{H^2},
\]
\[
B(\psi_1, \psi_1)|= \|\psi_1\|_{H^2}^2+\|\psi_1\|_{H^\theta}^2\geq \|\psi_1\|_{H^2}^2
\]
with some constant $C > 0$. Meanwhile, we define a bounded linear functional $F$ in $H^2({\bf R}^n)$
\[
	~~~~~\left\langle F, \psi \right\rangle=(f+g, \psi)_{H^\theta},~~~~\forall~ \psi \in H^2({\bf R}^n)\subset H^\theta({\bf R}^n).
\]
By relying on the Lax-Milgram theorem, there exists a unique (weak) solution $u\in H^2({\bf R}^n)$ satisfying  
\[
~~~~~B(u,\psi)=\left\langle F, \psi \right\rangle,~~~~\forall~ \psi \in H^2({\bf R}^n)\subset H^\theta({\bf R}^n),
\]
which is consistent with the weak form of \eqref{d69}. Here we just proved the existence of the (weak) solution $u \in H^2({\bf R}^n)$ of \eqref{d69}. 

Finally, we have to show that $u \in H^{4-\theta}({\bf R}^n) = \mathcal{D}(P)$.
Taking  the Fourier transform on both sides of \eqref{d69}, we have
\begin{equation}\label{d27}
	\widehat{Pu}=\hat{g}+\hat{f}-\hat{u}.
\end{equation}
Calculating the $H^\theta$-norm of \eqref{d27} yields
\begin{equation}\label{}
	\int_{{\bf R}^n_\xi} (1+\delta |\xi|^{2\theta })|\widehat{Pu} |^2
	~d\xi =\int_{{\bf R}^n_\xi}
	(	1+|\xi|^{2\theta})
	|\hat{f}+\hat{g}-\hat{u}|^2 ~ d\xi,
\end{equation}
which implies
\[
\|Pu\|_{H^\theta }^2=\|f+g-u\|_{H^\theta}^{2} \leq C(\|f\|_{H^\theta}^2+\|g\|_{H^\theta}^2+\|u\|_{H^\theta}^2)
\]
with some constant $C > 0$. By Lemma \ref{d7}, we can obtain
\begin{equation}\label{key}
	\|u\|_{H^{4-\theta} }\leq C  \|Pu\|_{H^\theta},
\end{equation}
with some constant $C > 0$, that is, $u\in H^{4-\theta}({\bf R}^n)$ ($\theta \in [0,2]$). Then $v$ can be determined by $v=u-f\in H^{2}({\bf R}^n)$.
Therefore, we complete the proof of the surjectivity of the mapping $\mathcal{I}-\mathcal{A}$. \\

By the Lumer-Phillips theorem, we conclude that the operator $\mathcal{A}$ generates a strongly continuous contraction semigroup  on $\mathcal {H}$.
\hfill
$\Box$

Let us finalize the proof of Proposition \ref{prp1}.

\noindent{\it{Proof of Proposition \ref{prp1}.}}
By \eqref{d18}, Lemma \ref{d29} and semigroup theory of linear operators, we just need to prove that $\mathcal{B}U$ is bounded in $\mathcal{H}$. In fact, it follows from the definition of $\mathcal{B}$ that 
\begin{align}
 \Vert \mathcal{B} U \Vert^2_{H^{2}\times H^{\theta}}&=\Vert(I+\delta A^{\theta})^{-1}v \Vert_{H^{\theta}}^{2}\\
 &=\int_{{\bf R}^{\xi}_{n}}(1+\delta|\xi|^{2\theta})\left|  \frac{\hat{v}}{1+\delta|\xi|^{2\theta}}\right|^2 d\xi \notag \\
 &=\int_{{\bf R}^{\xi}_{n}}\frac{|\hat{v}|^2}{1+\delta|\xi|^{2\theta}}d\xi \notag \\
 & \leq
\int_{{\bf R}^{\xi}_{n}}{|\hat{v}|}^{2}d\xi \notag\\
&\leq
\Vert v \Vert_{H^{\theta}}^{2},
\end{align}
which completes the proof of the unique existence of mild solutions. The existence of so-called strong solutions is a consequence of standard semigroup theory. }
\hfill
$\Box$


\begin{thebibliography}{99}






\bibitem{Ba} J. Barrera and H. Volkmer, Asymptotic expansion of the $L^{2}$-norm of a solution of the strongly damped wave equation in space dimension $1$ and $2$, Asymptotic Anal. {\bf 121} (2021), no. 3-4, 367--399.

\bibitem{Be} W. Beckner, Pitt's inequality with sharp convolution estimates, Proc. Amer. Math. Soc. {\bf 136}, no. 5,  (2008), 1871--1885.  

\bibitem{B} P. Brenner, On $L_{p}-L_{p'}$ estimates for the wave equation, Math. Z. {\bf 145} (1975), 251--254.  

\bibitem{cr2013} R. C. Char\~ao, C. R. da Luz and R. Ikehata, New decay rates for a problem of plate dynamics with fractional damping, J Hyperbolic Differ Equ.  {\bf 10} (2013), 563--575.

\bibitem{rcc2020} R. C. Char\~ao and R. Ikehata, A note on decay rates of the local energy for wave equations with Lipschitz wavespeeds, J. Math. Anal. Appl. {\bf483} (2020) 123636.

\bibitem{CI} W. Chen and R. Ikehata, The Cauchy problem for the Moore-Gibson-Thompson equation in the dissipative case, J. Diff. Eqns {\bf 292} (2021), 176--219.


\bibitem{CT} W. Chen and H. Takeda, Large-time asymptotic behaviors for the classical thermoelastic system, J. Diff. Eqns {\bf 377} (2023), 809--848. 

\bibitem{DH} E. B. Davies and A. M. Hinz, Explicit constants for Rellich inequalities in $L^{p}(\Omega)$, Math. Z. {\bf 227} (1998), 511--523.

\bibitem{DS} W. Dan and Y. Shibata, On a local energy decay of solutions of a dissipative wave equation, Funkcial. Ekvac. {\bf 38} (1995), no. 3, 545--568.

\bibitem{FIM} T. Fukushima, R. Ikehata and H. Michihisa, Thresholds for low regularity solutions to wave equations with structural damping, J. Math. Anal. Appl. {\bf 494} (2021), 124669.

\bibitem{H2016} J. L. Horbach, Existencia De Soluc\~oes E Comportamento Assint\'otico \'Otimo Para Equac\~oes Dissipativas Tipo Placas/Boussinesq Generalizadas Em ${\bf R}^{n}$. Doctoral Thesis: Federal University of Santa Catarina (2016).

\bibitem{HIC} J. L. Horbach, R. Ikehata and R. C. Char\~ao, Optimal decay rates and asymptotic profile for the plate equation with structural damping. J. Math. Anal. Appl. {\bf 440} (2016), no. 2, 529--560. 

\bibitem{ikawa} M. Ikawa, Hyperbolic differential equations and wave phenomena; 2000. Translations of Mathematical Monographs, AMS.

\bibitem{I-04} R. Ikehata, New decay estimates for linear damped wave equations and its application to nonlinear problem, Math. Meth. Appl. Sci. {\bf 27} (2004), 865-889. doi: 10.1002/mma.476.


\bibitem{ike2005} R. Ikehata, Local energy decay for linear wave equations with variable coefficients, J. Math. Anal. Appl. {\bf306} (2005), 330--348.

\bibitem{I-14} R. Ikehata, Asymptotic profiles for wave equations with strong damping, J. Diff. Eqns {\bf 257} (2014), 2159-2177.

\bibitem{JHDE-ike} R. Ikehata, $L^{2}$-blowup estimates of the wave equation and its application to local energy decay, J. Hyperbolic Differ. Equ. {\bf 20}, no.1, (2023), 259--275. 

\bibitem{FE-ike} R. Ikehata, $L^2$-blowup estimates of the plate equation, Funk. Ekvac., in press (2023).

\bibitem{IM} R. Ikehata and T. Matsuyama, $L^{2}$-behaviour of solutions to the linear heat and wave equations in exterior domains, Sci. Math. Japon. {\bf 55} (2002), 33-42.

\bibitem{IO} R. Ikehata and M. Onodera, Remarks on large time behavior of the $L^{2}$-norm of solutions to strongly damped wave equations, Differ. Integral Equ. {\bf 30} (2017), 505--520.

\bibitem{ir2015} R. Ikehata and M. Soga, Asymptotic profiles for a strongly damped plate equation with lower order perturbation, Commun. Pure Appl. Anal. {\bf 14} (2015), 1759--1780.

\bibitem{IT} R. Ikehata and H. Takeda, Asymptotic profiles of solutions for structural damped wave equations, J. Dynam. Diff. Eqns {\bf 31} (2019), 537--571. 

\bibitem{IL} X. Li and R. Ikehata, $L^{2}$-growth property for wave equations with higher derivative terms, arXiv: 2307.13329v1 [math.AP] 25 Jul 2023.

\bibitem{lz1999} Z. Liu Z and S. Zheng, Semigroups Associated with Dissipative Systems, Chapman and Hall/CRC Research Notes in Mathematics, vol. 398. Boca
Raton, FL: Chapman and Hall/CRC; 1999.

\bibitem{cr2009} C. R. da Luz and R. C. Char\~ao, Asymptotic properties for a semilinear plate equation in unbounded domains, J Hyperbolic Differ Equ.  {\bf 6} (2009), 269--294.

\bibitem{MOW} S. Machihara, T. Ozawa and H. Wadade, Hardy type inequalities on balls, Tohoku Math. J. {\bf 65} (2013), 321--330.

\bibitem{mh2019} H. Michihisa, New asymptotic estimates of solutions for generalized Rosenau equations, Math. Methods Appl. Sci. {\bf 42} (2019), 4516--4542.

\bibitem{Mi} H. Michihisa, Optimal leading term of solutions to wave equations with strong damping terms, Hokkaido Math. J. {\bf 50} (2021), 165--186.

\bibitem{M} A. Miyachi, On some estimates for the wave equation in $L_p$ and $H_p$, J. Fac. Sci. Univ. Tokyo Sect. IA Math. {\bf 27} (1980), no. 2, 331--354.

\bibitem{p1999} M. A. Park, On some nonlinear nonlinear dispersive equations, Contemp. Math. {\bf 221} (1999) 211--216, AMS, Providence, RI, 1999.

\bibitem{p1999-2} M. A. Park, On the Rosenau equation, Mat. Aplic. Comp. {\bf 9} (1990), 145--152.

\bibitem{P} H. Pecher, $L^{p}$-Absch\"atzungen und klassische L\"osungen f\"ur nichtlineare Wellengleichungen. I, Math. Z. {\bf 150} (1976), 159--183.

\bibitem{Pl} J. C. Peral, $L^{p}$ estimates for the wave equation, J. Funct. Analysis {\bf 36} (1980), 114--145. 

\bibitem{pn2009} N. Polat and E. Abdulkadir, Existence and blow-up of solution of Cauchy problem for the generalized damped multidimensional Boussinesq equation, J. Math. Anal. Appl. {\bf 349} (2009), 10--20.

\bibitem{rp1968} P. Rosenau, A quasi-continuous description of a nonlinear transmission line, Physica Scripta {\bf34} (1986), 827--829.

\bibitem{rp1968-2} P. Rosenau, Dynamics of dense discrete systems. Prog Theor Phys. {\b 79} (1988), 1028--1042.

\bibitem{S} R. S. Strichartz, A priori estimates for the wave equation and some applications, J. Funct. Anal. {\bf 5} (1970), 519--531.

\bibitem{th2012} H. Takeda and S. Yoshikawa, On the initial value problem of the semilinear beam equation with weak damping II: Asymptotic profiles, J Differ.
Equ. {\bf253} (2012), 3061-3080.

\bibitem{W} W. Von Wahl, $L^{p}$-decay rates for homogeneous wave equations, Math. Z. {\bf 120} (1971), 93--106.

\bibitem{wangxu09} S. Wang and  G. Xu, The Cauchy problem for the Rosenau equation, Nonlinear Anal. Theory, Methods and Appl. {\bf 71} (2009), 456--466.

\bibitem{ws2013} H. Wang and S. Wang, Global existence and asymptotic behavior of solution for the Rosenau equation with hydrodynamical damped term, J. Math. Anal. Appl. {\bf 401} (2013), 763--773.

\end{thebibliography}
\end{document}